\documentstyle[11pt]{article}

\def\init{\setcounter{equation}{0}}

\newtheorem{theoreme}{Theorem}[section]
\newtheorem{proposition}[theoreme]{Proposition}
\newtheorem{lemme}[theoreme]{Lemma}

\def\rr{{\bf R}}

\def\zz{{\bf Z}}

\def\calt{{\cal T}}

\def\ie{{\it i.e., }}

\def\proof{{\it  Proof. }}
\def\qed{$\Box$}

\def\half{\frac{1}{2}}

\newcommand{\beq}{\begin{equation}}
\newcommand{\eeq}{\end{equation}}
\newcommand{\bet}{\begin{theoreme}}
\newcommand{\eet}{\end{theoreme}}
\newcommand{\bear}[1]{\begin{array}{#1}}
\newcommand{\ear}{\end{array}}

\def\aa{\alpha}
\def\l{\lambda}

\begin{document}

\title{Fuglede's conjecture for a union of two intervals}
\author{I. \L aba\\ Department of Mathematics\\Princeton University\\
Princeton, NJ 08544\\ U.S.A.\\{\it laba@math.princeton.edu}}
\maketitle

\begin{abstract}
We prove that a union of two intervals in $\rr$ is a spectral set if
and only if it tiles $\rr$ by translations. Mathematics Subject 
Classification: 42A99.
\end{abstract}

\section{The results}
\label{intro}
\init

A Borel set $\Omega\subset\rr^n$ of positive measure is said to 
{\it tile $\rr^n$ by translations} if there is a discrete set 
$T\subset\rr^n$ such that, up to sets of measure 0, the sets $\Omega
+t,\ t\in T,$ are disjoint and $\bigcup_{t\in T}(\Omega+t)=\rr^n$.
We may rescale $\Omega$ so that $|\Omega|=1$.
We say that $\Lambda=\{\lambda_k:\ k\in\zz\}\subset\rr^n$ is
a {\it spectrum} for $\Omega$ if:
\beq
\{e^{2\pi i\l_k\cdot x}\}_{k\in\zz}\hbox{ is an orthonormal basis for }
L^2(\Omega).
\label{a.e00}
\eeq
A {\it spectral set} is a domain $\Omega\in\rr^n$ such that (\ref{a.e00})
holds for some $\Lambda$.  

Fuglede \cite{Fug} conjectured that a domain
$\Omega\subset\rr^n$ is a spectral set if and only if it tiles $\rr^n$
by translations, and proved this conjecture under the assumption that
either $\Lambda$ or $T$ is a lattice.  The conjecture is related
to the question of the existence of commuting self-adjoint extensions
of the operators $-i\frac{\partial}{\partial x_j}$, $j=1,\dots,n$ 
\cite{Fug}, \cite{J}, \cite{P1}; other relations between the tiling
and spectral properties of subsets of $\rr^n$ have been conjectured
and in some cases proved, see \cite{IP}, \cite{JP1}, \cite{JP2},
\cite{K2}, \cite{LRW}, \cite{LW2}.

Recently there has been significant progress on the special case of
the conjecture when $\Omega$ is assumed to be convex \cite{K1},
\cite{IKP}, \cite{IKT1}, and in particular the 2-dimensional convex case
appears to be nearly resolved \cite{IKT2}. The non-convex case is
considerably more complicated and is not understood even in dimension 1.
The strongest results yet in that direction seem to be those of Lagarias
and Wang \cite{LW1}, \cite{LW2}, who proved that all tilings
of $\rr$ by a bounded region must be periodic and that the corresponding
translation sets are rational up to affine transformations, which in turn
leads to a structure theorem for bounded tiles.  It was also observed in
\cite{LW2} that the ``tiling implies spectrum" part of Fuglede's
conjecture for compact sets in $\rr$ would follow from a conjecture
of Tijdeman \cite{Tij} concerning factorization of finite cyclic groups;
however, Tijdeman's conjecture is now known to fail without additional
assumptions \cite{CM}. 
See also \cite{New}, \cite{CM} for partial results on the related problem
of characterizing all tilings of $\zz$ by a finite set, and \cite{LW2},
\cite{P2} for a classification of domains in $\rr^n$ which have
$L+\zz^n$ as a spectrum for some finite set $L$.

The purpose of the present article is to address the following special
case of Fuglede's conjecture in one dimension.
Let $\Omega=I_1\cup I_2$, where $I_1,I_2$ are disjoint intervals
of non-zero length.  By scaling, translation, and symmetric 
reflection, we may assume that:
\beq
\Omega=(0,r)\cup(a,a+1-r),\ 0< r\leq \half,\  a\geq r.
\label{e.omega1}
\eeq
Our first theorem characterizes all $\Omega$'s of the form 
(\ref{e.omega1}) which are spectral sets.

\begin{theoreme}
Suppose that $\Lambda$ is a spectrum for $\Omega$, $0\in\Lambda$.
Then at least one of the following holds:

\smallskip
(i) $a-r\in\zz$ and $\Lambda=\zz$;

\smallskip
(ii) $r=\half$, $a=\frac{n}{2}$ for some $n\in\zz$, and 
$\Lambda=2\zz\bigcup(\frac{p}{n}+2\zz)$ for some odd integer $p$.

\smallskip\noindent
Conversely, if $\Omega$, $\Lambda$ satisfy (\ref{e.omega1}) and if
either {\it(i)} or {\it(ii)} holds, then $\Lambda$ is a spectrum 
for $\Omega$.
\label{s.thm1}
\end{theoreme}

As a corollary, we prove that Fuglede's conjecture is true for a
union of two intervals.

\begin{theoreme}
Let $\Omega\subset\rr$ be a union of two disjoint intervals,
$|\Omega|=1$.  Then $\Omega$ has a spectrum if and only if it tiles
$\rr$ by translations.
\label{s.thm2}
\end{theoreme}

Theorem \ref{s.thm2} follows easily from Theorem \ref{s.thm1}.
We may assume that $\Omega$ is as in (\ref{e.omega1}).
Suppose that $\Lambda$ is a spectrum for $\Omega$; without loss of
generality we may assume that $0\in\Lambda$.  Then by Theorem 
\ref{s.thm1} one of the conclusions {\it(i), (ii)} must hold,
and in each of these cases $\Omega$ tiles $\rr$ by translations.
Conversely, if $\Omega$ tiles $\rr$ by translations, by Proposition
\ref{tis.prop1} $\Omega$ must satisfy Theorem \ref{s.thm1}{\it(i)} 
or {\it(ii)}; the second part of Theorem \ref{s.thm1} implies then 
that $\Omega$ has a spectrum.

Theorem \ref{s.thm1} will be proved as follows. Suppose that $\Lambda
=\{\l_k:\ k\in\zz\}$is a spectrum for $\Omega$; we may assume that 
$\l_0=0$. Let $\l_{kk'}=\l_k
-\l_{k'}$, $\Lambda-\Lambda=\{\l_{kk'}:\ k,k'\in\zz\}$, and:
\beq
Z_\Omega=\{0\}\cup\{\l\in\rr:\ \hat\chi_\Omega(\lambda)=0\}.
\label{e.Z}
\eeq
Then the functions $e^{2\pi i\l_k x}$ are mutually orthogonal in 
$L^2(\Omega)$, hence $\Lambda\subset\Lambda-\Lambda\subset Z_\Omega$. 
This will lead to a number of restrictions on the possible values
of $\lambda_k$.  Next, let:
\beq
\phi_\l(x)=\chi_{(0,r)}e^{2\pi i\l x},
\label{s.e11}
\eeq
where $\chi_{(0,r)}$ denotes the characteristic function of $(0,r)$. 
By Parseval's formula, the Fourier coefficients 
$c_k=\int_0^r e^{2\pi i(\l-\l_k)x}dx$ of $\phi_\l$ satisfy:
\beq
\sum_{k\in\zz}c_k^2=\|\chi_{(0,r)}e^{2\pi i\l x}\|_{L^2(\Omega)}^2=r.
\label{s.e10}
\eeq
Given that the $\lambda_k$'s are subject to the orthogonality
restrictions mentioned above, we will find that there are not
enough $\lambda_k$'s for (\ref{s.e10}) to hold unless the conditions
of Theorem \ref{s.thm1} are satisfied.

The author is grateful to Alex Iosevich for helpful conversations 
about spectral sets and Fuglede's conjecture.

\section{Tiling implies spectrum}
\label{tis}
\init

\begin{proposition}
If $\Omega$ as in (\ref{e.omega1}) tiles $\rr$ by translations, it
must satisfy {\it(i)} or {\it(ii)} of Theorem \ref{s.thm1}.
\label{tis.prop1}
\end{proposition}

\proof
Suppose that $\rr$ may be tiled by translates of $\Omega$. Assume 
first that $r=\half$. Any copy of $\Omega$ used in the tiling has a
``gap" of length $a-r=a-\half$, which must be covered by non-overlapping
intervals of length $\half$; hence $a\in\half\zz$ as in Theorem 
\ref{s.thm1}{\it(ii)}.

Assume now that $0<r<\half$. Let $I_1=(0,r)$, $I_2=(a,a+1-r)$. 
We will prove that translates of $I_1$ and $I_2$ must alternate in any
tiling $\calt$ of $\rr$ by translates of $\Omega$; this implies
immediately that $a-r\in\zz$ as in Theorem \ref{s.thm1}{\it(i)}. 

\begin{itemize}

\item If $\calt$ contained two consecutive translates $(\tau,
\tau+r)$ and $(\tau+r,\tau+2r)$ of $I_1$, it would also 
contain the matching translates $(\tau+a, \tau+a+1-r)$ and 
$(\tau+a+r,\tau+a+1)$ of $I_2$, which is impossible since the
latter two intervals overlap. 

\item Suppose now that $\calt$ contains
two consecutive translates $(\tau+a,\tau+a+1-r)$ and $(\tau+a+1-r,
\tau+a+2-2r)$ of $I_2$; then $\calt$ must also contain the matching
translates $I'_1=(\tau,\tau+r)$ and $I''_1=(\tau+1-r,\tau+2-2r)$
of $I_1$. The gap between $I'_1$ and $I''_1$ has length $1-2r$,
which is strictly less than $1-r=|I_2|$, so that $I'_1$ must be 
followed by another translate of $I_1$. But this has just been 
shown to be impossible.
\qed

\end{itemize}

\bigskip

Next, we prove the second part of Theorem \ref{s.thm1}. This easy result
appears to have been known to several authors, see e.g., the examples
in \cite{Fug}, \cite{JP1}, \cite{LW2}. Since we will rely on it later on
in the proof of the ``hard" part of the theorem, we include the short 
proof. 

\begin{proposition}
If $\Lambda$ and $\Omega$ are as in Theorem \ref{s.thm1}{\it(i)}
or {\it(ii)}, then $\Lambda$ is a spectrum for $\Omega$. 
\label{tis.prop2}
\end{proposition}

\proof
If {\it(i)} holds, then $\Omega$ is a fundamental domain for $\zz$
and consequently $\Lambda=\zz$ is a spectrum \cite{Fug}.
Suppose now that {\it(ii)} holds. For any function $f$ on $\Omega$,
we define functions $f_+,f_-$:
\[
f_+(x)=\half(f(x)+f(x')),\ f_-(x)=\half(f(x)-f(x')),\ x\in \Omega,
\]
where $x'=x+a$ if $x\in(0,\half)$, and $x'=x-a$ if $x\in(a,a+\half)$.
Then:
\[
f(x)=f_+(x)+f_-(x),\ f_+(x)=f_+(x'),\ f_-(x)=-f_-(x').
\]
It therefore suffices to prove that:
\beq
g(x)=\sum_{k\in\zz}c_ke^{4k\pi ix}\hbox{ for any }g(x)
\hbox{ such that }g(x)=g(x'),
\label{s.e20}
\eeq
\beq
h(x)=\sum_{k\in\zz}c'_ke^{(4k+\frac{2p}{n})\pi ix}\hbox{ for any }
h(x)\hbox{ such that }h(x)=-h(x').
\label{s.e21}
\eeq
Since $e^{4k\pi ix}$, $k\in\zz$, is a spectrum for $(0,\half)$, we
have:
\[
g(x)=\sum_{k\in\zz}c_ke^{4k\pi ix},\ 
h(x)=e^{\frac{2p}{n}\pi ix}\sum_{k\in\zz}c'_ke^{4k\pi ix},
\ x\in(0,\half).
\]
(\ref{s.e20}) follows immediately by periodicity.
From the second equation above we find that 
(\ref{s.e21}) holds for all $x\in(0,\half)$, and that for such $x$:
\[
e^{\frac{2p}{n}\pi i(x+a)}\sum_{k\in\zz}c'_ke^{4k\pi i(x+a)}
=-e^{\frac{2p}{n}\pi ix}\sum_{k\in\zz}c'_ke^{4k\pi ix}
=-h(x)=h(x+a),
\]
where we used that $\frac{2p}{n}a=p$ is odd. Hence (\ref{s.e21})
holds also for $x\in(a,a+\half)$.
This ends the proof of Proposition \ref{tis.prop2}.
\qed

\section{Orthogonality}
\label{o}
\init

We now begin the proof of the first part of Theorem \ref{s.thm1}.
Throughout the rest of the paper, $\Omega$ is assumed to satisfy
(\ref{e.omega1}), $\Lambda=\{\l_k:\ k\in\zz\}$ is a spectrum for
$\Omega$, $\l_0=0$, $\l_{kk'}=\l_k-\l_{k'}$, $\Lambda-\Lambda=
\{\l_{kk'}:\ k,k'\in\zz\}$, and $Z_\Omega$ is
defined by (\ref{e.Z}).

\begin{lemme}
$Z_\Omega= Z_1\cup Z_2\cup Z_3$, where:
\[
\begin{array}{l}
Z_1=\{\l\in\rr:\ \l a\in \zz+\half,\ \l(2r-1)\in\zz\},
\\[3mm]
Z_2=\{\l\in\zz:\ \l r\in\zz\},
\\[3mm]
Z_3=\{\l\in\zz:\ \l(a-r)\in\zz\}.
\end{array}
\]
\label{s.lemma1}
\end{lemme}

\proof
Suppose that $\l\neq 0$, $\l\in Z_\Omega$. Then:
\[
\int_\Omega e^{2\pi i\l x}dx=
e^{2\pi i\l r}-1+e^{2\pi i\l(a+1-r)}-e^{2\pi i\l a}=0.
\]
All solutions to $z_1+z_2+z_3+1=0$, $|z_i|=1$, must be of the form 
$\{z_1,z_2,z_3\}=\{-1,z_*,-z_*\}$. Hence $\l\in Z_\Omega$ if and only
if one of the following holds.

\begin{itemize}

\item $e^{2\pi i\l a}=-1$ and $e^{2\pi i\l r}+e^{2\pi i\l(a+1-r)}=0$,
hence $\l\in Z_1$;

\item $e^{2\pi i\l r}=1$ and $e^{2\pi i\l(1-r)}=1$, hence $\l\in Z_2$;

\item $e^{2\pi i\l(a+1-r)}=1$ and $e^{2\pi i\l a}=e^{2\pi i\l r}$, 
hence $\l\in Z_3$. $\Box$

\end{itemize}

Observe that $Z_2$, $Z_3$ are additive subgroups of $\zz$.

\begin{lemme}
At least one of the following holds:
\beq
\Lambda\subset Z_1\cup Z_2,
\label{s.e1}
\eeq
\beq
\Lambda\subset Z_1\cup Z_3.
\label{s.e2}
\eeq
\label{s.lemma2}
\end{lemme}

\proof
By Lemma \ref{s.lemma1}, $\Lambda\subset\Lambda-\Lambda\subset Z_\Omega
\subset Z_1\cup Z_2\cup Z_3$.
If $Z_2\subset Z_3$, (\ref{s.e2}) holds; suppose therefore that
there is a $\lambda_i\in Z_2\setminus Z_3$. It suffices to prove
that for any $\lambda_j\in Z_3$ we must have $\l_j\in Z_1$ or
$\l_j\in Z_2$.

Let $\l_j\in Z_3$, then $\l_{ij}=\l_i-\l_j\in Z_\Omega$ by
orthogonality. By Lemma \ref{s.lemma1}, $\l_{ij}\in Z_1\cup
Z_2\cup Z_3$. If $\l_{ij}\in Z_2$, then $\l_j\in Z_2$ and we are
done, and if $\l_{ij}\in Z_3$, then $\l_i\in Z_3$, which 
contradicts our assumption. Assume therefore that $\l_{ij}
\in Z_1$. Then:
\[
\l_{ij}\in\zz,\ \l_{ij}a\in\zz+\half,
\ \l_{ij}(2r-1)\in\zz,
\]
hence:
\[
2\l_j r=2\l_i r-\l_{ij}(2r-1)-\l_{ij}\in\zz.
\]
If $\l_jr\in\zz$, then $\l_j\in Z_2$; if $\l_j r\in\zz+\half$,
then $\l_j a\in\zz+\half$ by the definition of $Z_3$ and 
$\l_j(2r-1)\in\zz$, so that $\l_j\in Z_1$.
$\Box$

\begin{lemme}
(i) $\Lambda\subset Z_2$ is not possible;

(ii) $\Lambda\subset Z_3$ is possible only if $a-r\in\zz$
and $\Lambda=Z_3=\zz$.
\label{s.lemma3}
\end{lemme}

\proof
Suppose that $\Lambda\subset Z_i$ for $i=2$ or $3$. Since $Z_i$
is an additive subgroup of $\zz$, we must have $Z_i=p\zz$ for
some integer $p>0$.  Furthermore, if there was a $\lambda\in p\zz
\setminus\Lambda$, we would have $\l_k-\l\in p\zz$ and hence 
$e^{2\pi i\l x}$ would be orthogonal to $e^{2\pi i\l_k x}$ for
all $\l_k\in\Lambda$, which would contradict (\ref{a.e00}).
Hence $\Lambda=Z_i=p\zz$.  We also observe that if $p$ was
$\geq 2$, any function of the form $f(x)=\sum_{k\in\zz}
c_ke^{2\pi i\l_k x}$ would be periodic with period $\frac{1}{p}
\leq\half$, which again would contradict (\ref{a.e00}).
Thus $\Lambda=Z_i=\zz$.

If $i=2$, this is not possible, since $nr$ cannot be integer for
all $n\in\zz$ if $r\leq\half$. If $i=3$, we obtain that $n(a-r)
\in\zz$ for all $n\in\zz$; letting $n=1$, we find that $a-r\in\zz$.
\qed

\bigskip

If $\Omega,\Lambda$ are as in Lemma \ref{s.lemma3}{\it(ii)}, then
Theorem \ref{s.thm1}{\it(i)} is satisfied and we are done.
Thus we may assume throughout the sequel that:
\beq
\Lambda\not\subset Z_2,\ \Lambda\not\subset Z_3.
\label{s.e5}
\eeq

\begin{lemme}
If (\ref{s.e5}) holds, then $\Lambda\subset Z_1\cup(Z_2\cap Z_3)$.
\label{s.lemma4}
\end{lemme}

\proof
By Lemma \ref{s.lemma2}), it suffices to prove that:
\beq
\hbox{if }\Lambda\cap(Z_1\setminus Z_2)\neq\emptyset,
\hbox{ then }\Lambda\cap Z_2\subset \Lambda\cap Z_3;
\label{s.e6}
\eeq
\beq
\hbox{if }\Lambda\cap(Z_1\setminus Z_3)\neq\emptyset,
\hbox{ then }\Lambda\cap Z_3\subset \Lambda\cap Z_2;
\label{s.e7}
\eeq
We will prove (\ref{s.e6}); the proof of (\ref{s.e7}) is almost
identical. Suppose that $\l_i\in Z_1\setminus Z_2$, and let
$\l_j\in Z_2$. By Lemma \ref{s.lemma1}, $\l_{ij}$ belongs to at
least one of $Z_1$, $Z_2$, $Z_3$; moreover, $\l_{ij}\in Z_2$
would imply $\l_i\in Z_2$ and contradict the above supposition.
Thus we only need consider the following two cases.

\begin{itemize}

\item Let $\l_{ij}\in Z_1$. Then $\l_i a,\l_{ij}a\in
\zz+\half$, hence $\l_j a\in\zz$ and $\l_j\in Z_2\cap Z_3$.

\item Assume now that $\l_{ij}\in Z_3$. Then $\l_i\in\zz$, 
hence $2\l_i r\in\zz$. We cannot have $\l_i r\in \zz$, since
then $\l_i$ would be in $Z_2$; therefore $\l_i r\in\zz+\half$.
Hence $\l_i(a-r)\in\zz$; since also $\l_{ij}(a-r)\in\zz$, we
obtain that $\l_j(a-r)\in\zz$ and $\l_j\in Z_2\cap Z_3$.
$\Box$.

\end{itemize}

\begin{lemme} Assume (\ref{s.e5}). Then:

\smallskip
(i) $\Lambda-\Lambda\subset Z_1\cup(Z_2\cap Z_3)$;

\smallskip
(ii) $\Lambda\cap Z_1\subset\l_*+r^{-1}\zz$ for some $\l_*\in\rr$.
\label{s.lemma5}
\end{lemme}

\proof
For $k\in\zz$, let $\Lambda_k=\Lambda-\l_k=\{\l_{jk}:\ j\in\zz\}$. 
Then $\Lambda_k$ is also a spectrum for $\Omega$ and $0\in
\Lambda_k$, hence all of the results obtained so far apply with
$\Lambda$ replaced by $\Lambda_k$. Thus {\it(i)} follows from
Lemmas \ref{s.lemma3} and \ref{s.lemma4}.

To prove {\it(ii)}, it suffices to verify that $\l_{ij}r\in\zz$
whenever $\l_i,\l_j\in Z_1$. Indeed, if $\l_i,\l_j\in Z_1$, then
$\l_{ij}a\in\zz$, hence $\l_{ij}\notin Z_1$ and therefore,
by {\it(i)}, $\l_{ij}\in Z_2\cap Z_3$. But this implies that
$\l_{ij}r\in\zz$.
\qed

\section{Completeness}
\label{c}
\init

Fix $j,n\in\zz$, and consider the function $\phi_\lambda$ defined
by (\ref{s.e11}) with $\lambda=\l_j-nr^{-1}$. The Fourier coefficients
of $\phi_{\lambda}$ are:
\[
c_k=\int_0^r e^{2\pi i(\l-\l_k)x}dx
=\int_0^r e^{2\pi i(\l_{jk}-nr^{-1})x}dx,
\]
hence $c_k=r$ if $\l_{jk}=nr^{-1}$, and:
\beq
c_k=\frac{1}{2\pi i(\l_{jk}-nr^{-1})}
\Big(e^{2\pi i(\l_{jk}r-n)}-1\Big),\ \l_{jk}\neq nr^{-1}.
\label{c.e1}
\eeq
Define $\aa_{jk}=\l_{jk}r$. Plugging (\ref{c.e1}) into 
(\ref{s.e10}), we obtain that for all $j\in\zz$:
\begin{equation}
\frac{1}{r}=1+\sum_{k:\aa_{jk}\notin\zz}
\frac{1}{4\pi^2\alpha_{jk}^2} \Big|e^{2\pi i\aa_{jk}}-1\Big|^2,
\label{c.e2}
\end{equation}
and for all $n,j\in\zz$:
\begin{equation}
\frac{1}{r}=\delta_{n,j}+\sum_{k:\aa_{jk}\notin\zz}\frac{1}{4\pi^2
(\alpha_{jk}-n)^2}\Big|e^{2\pi i(\aa_{jk}-n)}-1\Big|^2,
\label{c.e3}
\end{equation}
where $\delta_{n,j}=1$ if there is a $k\in\zz$ such that 
$\aa_{jk}=n$, and $\delta_{n,j}=0$ otherwise.

We define the equivalence relation between the indices $k,k'$:
\[
k\sim k' \ \Leftrightarrow\ \ \aa_{kk'}\in\zz,
\]
and denote by $A_1,A_2,\dots,A_m,\dots\subset\zz$ the (non-empty and
disjoint) equivalence classes with respect to this relation.
Hence $k,k'$ belong to the same $A_m$ if and only if $\aa_{kk'}\in
\zz$; in particular, $A_m\subset \beta_m+\zz$ for some $\beta_m
\in[0,1).$

\begin{lemme}
Let $M$ denote the number of distinct and non-empty $A_m$'s. Then:
\beq
M\geq r^{-1}.
\label{c.e4}
\eeq
Moreover, if one of the $A_m$'s skips a number (\ie $A_m\neq
\beta_m+\zz$), then $M\geq r^{-1}+1$.
\label{c.lemma1}
\end{lemme}

\bigskip

\proof
For each $m,m'$, let $\beta_{mm'}=\beta_m-\beta_{m'}$; note
that $\beta_{mm'}\neq 0$ if $m\neq m'$. Fix $m'$ and $j\in A_{m'}$,
then (\ref{c.e2}) may be rewritten as:
\beq
\frac{1}{r}=1+\sum_{m\neq m'}S_{mm'},
\label{c.e5}
\eeq
where:
\[
S_{mm'}=\sum_{k\in A_m}\frac{1}{4\pi^2\aa_{jk}^2}
\Big|e^{2\pi i\beta_{mm'}}-1\Big|^2.
\]
Clearly:
\beq
S_{mm'}\leq \tilde S(\beta_{mm'}),
\label{c.e7}
\eeq
where:
\beq
\tilde S(\beta)=\sum_{k\in \zz}\frac{1}{4\pi^2(\beta+k)^2}
\Big|e^{2\pi i\beta}-1\Big|^2.
\label{c.e6}
\eeq
Hence (\ref{c.e4}) follows from (\ref{c.e5}) and Lemma \ref{c.lemma2}
below. 

Suppose now that $A_{m'}$ skips a number. Then we may find $j\in
A_{m'}$ and $n\in\zz$ such that $\delta_{n,j}=0$, and (\ref{c.e4}) 
may be improved to $M\geq 1+r^{-1}$ by using (\ref{c.e3}) instead
of (\ref{c.e2}). 
\qed

\begin{lemme}
Let $\tilde S(\beta)$ be as in (\ref{c.e6}), then $\tilde S(\beta)
=1$ for all $0<\beta<1$.
\label{c.lemma2}
\end{lemme}

\proof
By Proposition \ref{tis.prop2}, $\Lambda=2\zz\cup(\frac{p}{n}+2\zz)$,
where $n\in\zz$ and $p$ is an odd integer, is a spectrum for 
$\Omega=(0,\half)\cup(\frac{n}{2},\frac{n+1}{2})$. Plugging this back
into (\ref{c.e2}) we obtain that: 
\[
1=\sum_{k\in\zz}\frac{1}{4\pi^2(\beta+k)^2}
\Big|e^{2\pi i\beta}-1\Big|^2
\]
for $\beta=\frac{p}{2n}$.  However, the set of $\beta$ of this form 
is dense in $\rr$, hence by continuity the lemma holds for all
$\beta\in (0,1)$.
\qed

\section{Conclusion}
\label{p}
\init

{\it Proof of Theorem \ref{s.thm1}.}
If $\Lambda$ is as in Lemma \ref{s.lemma3}{\it(ii)}, then
Theorem \ref{s.thm1}{\it(i)} is satisfied; we may therefore 
assume that (\ref{s.e5}) holds.  From Lemma \ref{s.lemma5} we have:
\beq
\Lambda-\Lambda\subset Z_1\cup(Z_2\cap Z_3),\ 
Z_2\cap Z_3\subset r^{-1}\zz,\ Z_1\subset(\l_*+r^{-1}\zz),
\label{s.e15}
\eeq
for some $\l_*\in\rr$, hence $M\leq 2$. However, by Lemma \ref{c.lemma1}
$M\geq r^{-1}\geq 2$, and this may be improved to $M\geq 3$ if one
of the $A_m$'s skips a number. Therefore we must have $r=\half$ and:
\beq
\Lambda-\Lambda=2\zz\cup(\l_*+2\zz),\ 
Z_2\cap Z_3=2\zz,\ Z_1=\l_*+2\zz.
\label{s.e16}
\eeq
Pick $\l_{ij},\l_{kl}\in Z_1$ such that $\l_{ij}-\l_{kl}=2$.
From the definition of $Z_1$ we have
$\l_{ij}a,\l_{kl}a\in\zz+\half$, hence:
\[
2a=\frac{a}{r}=\l_{ij}a-\l_{kl}a\in\zz,
\]
so that $a=\frac{n}{2}$ for some $n\in\zz$. Finally, we have
$\l_*a=\half n\l_*\in\zz+\half$, hence $\l_*n=p$ for some
odd integer $p$. Thus $\Omega$ and $\Lambda$ satisfy {\it(ii)}
of Theorem \ref{s.thm1}.


\end{document}